\documentclass{ifacconf}
\usepackage{graphicx}      

\usepackage{graphics} 
\usepackage{epsfig} 
\usepackage{mathptmx} 
\usepackage{times} 
\usepackage{amsmath} 
\usepackage{amssymb}  
\usepackage{algorithm}
\usepackage[usenames, dvipsnames]{color}
\usepackage{algpseudocode}
\usepackage{natbib}        
\usepackage{psfrag}
\usepackage{diagbox}

\usepackage{subcaption}
\newlength\figwidth
\usepackage[
  tableposition = top, 
  width = 0.48\textwidth 
]{caption}

\begin{document}
\begin{frontmatter}

\title{A Penalty Method Based Approach for Autonomous Navigation using Nonlinear Model Predictive Control\thanksref{footnoteinfo}} 

\thanks[footnoteinfo]{This work benefits from KU Leuven-BOF PFV/10/002 Centre of Excellence: Optimization in Engineering (OPTEC), from the project G0C4515N of the Research Foundation-Flanders (FWO-Flanders), from Flanders Make ICON project: Avoidance of collisions and obstacles in narrow lanes, and from the KU Leuven Research project C14/15/067: B-spline based certificates of positivity with applications in engineering, from KU Leuven internal funding: StG/15/043, from Fonds de la Recherche Scientifique -- FNRS and the Fonds Wetenschappelijk Onderzoek -- Vlaanderen under EOS Project no 30468160 (SeLMA) and from FWO projects G086318N, G086518N.}

\author[First]{Ben Hermans} 
\author[Second]{Panagiotis Patrinos} 
\author[First]{Goele Pipeleers}
\address[First]{MECO Research Team, Department of Mechanical Engineering, KU Leuven, Leuven, Belgium \\ DMMS lab, Flanders Make, Leuven, Belgium \\ 
(e-mail: ben.hermans2@kuleuven.be).}
\address[Second]{Division ESAT-STADIUS, Department of Electrical Engineering, KU Leuven, Leuven, Belgium}

\begin{abstract}                
This paper presents a novel model predictive control strategy for  controlling autonomous motion systems moving through an environment with obstacles of general shape. In order to solve such a generic non-convex optimization problem and find a feasible trajectory that reaches the destination, the approach employs a quadratic penalty method to enforce the obstacle avoidance constraints, and several heuristics to bypass local minima behind an obstacle. The quadratic penalty method itself aids in avoiding such local minima by gradually finding a path around the obstacle as the penalty factors are successively increased. The inner optimization problems are solved in real time using the proximal averaged Newton-type method for optimal control (PANOC), a first-order method which exhibits low runtime and is suited for embedded applications. The method is validated by extensive numerical simulations and shown to outperform state-of-the-art solvers in runtime and robustness. 

\end{abstract}


\end{frontmatter}

\section{Introduction}


Driverless cars, fruit-picking robots and automated guided vehicles in a warehouse are examples of autonomous motion systems that are up-and-coming in industry. In all these applications, the computation of a collision-free trajectory is essential. Computing motion trajectories that satisfy collision-avoidance constraints has been the topic of substantial research, resulting in a variety of methods, including the potential field method \citep{ge2002dynamic, montiel2015path} and methods using velocity obstacles \citep{guy2009clearpath}. Another popular approach relies on constructing a graph by discretizing the geometric domain, using vonoroi diagrams \citep{takahashi1989motion} or a simple grid of square elements, and performing a subsequent graph search. These graph searches are typically performed using Dijkstra's algorithm or one of its variants with additional heuristics, such as the popular A* \citep{hart1968formal}. Grid based motion planning problems can also be solved by wavefront planners, such as D* \citep{stentz1994optimal} and its variants. The main drawback of graph search methods is that they usually do not consider kinematic constraints of the motion vehicle, which is a problem especially for nonholomic vehicles. Recently, optimization-based strategies for solving motion planning problems, such as Model Predictive Control (MPC), are also becoming more popular.

MPC is a control strategy in which an optimal control problem is solved at every time instant \citep{rawlings2009model}. Of the resulting optimal control sequence, the first one is applied to the plant, and the procedure is repeated. A significant advantage of MPC is its ability to take into account constraints on the inputs and states, such as collision-avoidance constraints.

One of the main challenges of the practical application of MPC is the strict real-time constraint for solving the optimal control problems. High sampling frequencies are typically required for the system to be able respond appropriately to disturbances and changes in the environment. Moreover, solvers often have to run on resource-constrained embedded hardware. The traditional solvers for numerical optimization, Sequential Quadratic Programming (SQP) and Interior Point (IP) methods, are not very suitable for this purpose as they require the costly operation of solving a linear system of equations at every iteration. In contrast, first-order methods do not require this operation and often involve only simple steps. This explains their increasing popularity for solving MPC problems, \citep{richter2012computational,patrinos2014accelerated,jerez2014embedded}.

Furthermore, the optimization algorithm and the problem form, in particular the constraints, are generally linked. For example, most SQP solvers assume that the linear independence constraint qualification (LICQ) is satisfied. 

Often only simple obstacle shapes, such as circular \citep{wang2014synthesis} and rectangular obstacles are considered. Another approach is based on the separating hyperplane theorem \citep{boyd2004convex}, and allows for the separation of a convex motion system and convex obstacles, or between convex motion systems, as illustrated by \citep{debrouwere2013time} and \citep{mercy2017spline}. Recently, \cite{embedded} have proposed a novel constraint formulation to incorporate general obstacle shapes, described as the intersection of a set of nonlinear inequalities, in the optimization problem.

This paper embeds the obstacle constraint formulation presented by \citet{embedded} in a penalty method framework to calculate a trajectory while satisfying collision-avoidance constraints. The penalty parameters allow for  a trade-off between the optimality of the trajectory and the extent to which the obstacle constraints may be violated. Virtual enlargements ensure that this trade-off results in a trajectory that avoids all real obstacles. Moreover, the application of the penalty method lowers the likelihood of getting stuck in local optima due to obstacles, as the trajectory is gradually formed around the obstacle. In addition, some heuristics are developed for dealing further with these local optima.

The resulting optimization problems are solved using the proximal averaged Newton-type method for optimal control (PANOC), as proposed in \citep{stella2017simple}. As this method combines projected gradient and limited-memory quasi-Newton steps, its implementation is simple and it can achieve a fast rate of convergence. An automatic differentiation toolbox, CasADi \citep{Andersson2018}, is used to efficiently compute the value of the objective function and its derivative.

The proposed algorithm is validated for a set of obstacle configurations and for different vehicle models. It is shown to be successful in avoiding obstacles of arbitrary shape, as long as they can be described by the intersection of a set of nonlinear inequalities. In addition, our algorithm is benchmarked against state-of-the-art SQP and IP methods, and is found to outperform these methods both in terms of runtime and robustness. 

This paper is organized as follows: Section II describes the obstacle constraint formulation introduced by \cite{embedded}. In addition, the resulting optimization problem is presented. Section III discusses the methodology for solving this problem, consisting of necessary reformulations, the first-order algorithm PANOC, the quadratic penalty method and several heuristics. Section IV shows and discusses the  numerical simulations results. Section V draws the main conclusions of the paper.

\section{PROBLEM FORMULATION}
An obstacle constraint formulation that can deal with general obstacle shapes was first introduced by \cite{embedded}. This constraint formulation is an essential part of the optimization problem considered in this paper. Another element of the problem are the kinematics of the motion system for which trajectories are calculated. Both these elements are analyzed in the subsections below and finally incorporated in a nonlinear model predictive control (NMPC) framework. 
\subsection{Obstacle constraint formulation}
In this work, we consider obstacle shapes that can be defined as the intersection of a set of $m$ nonlinear inequalities:
\begin{equation}
O = \{z \in {\rm I\!R}^{n_z} : h_i(z)>0, \hspace{0.1cm} i=1,...,m\}.
\end{equation}
Here, $z$ denotes the position vector, $n_z$ the number of dimensions, and $h_i : {\rm I\!R}^{n_z} \rightarrow {\rm I\!R}$ are continuously differentiable functions with Lipschitz continuous gradients. In the remainder of this paper the considered geometry will be two dimensional, thus $n_z = 2$ and $z = (x,y)$. To take into account moving obstacles, time dependent functions $h_i(z,t)$ can also be used in this formulation. However, this paper will only consider the static case. The obstacle avoidance constraint ($z \not\in O$) can then be written as follows:
\begin{equation}
\exists i \in {1,...,m} : \hspace{0.5cm} h_i(z) \leq 0.
\end{equation}
In other words, at least one of the inequalities defining the obstacle  must be violated. This condition can be rewritten as the following equality constraint:
\begin{equation} \label{eq:Obstcost}
\psi(z):=\displaystyle\prod ^{m}_{i=1} [h_i(z)]_+ = 0,
\end{equation}
where the operator $[h_i(z)]_+$ is defined as $\max(h_i(z),0)$. 

The obstacle avoidance constraint (\ref{eq:Obstcost}) is linked to vertical complementary constraints \citep{scheel2000mathematical}, 
as it can be rewritten as 
\begin{equation}
\mathrm{min}([h_1(z)]_+^2,..., [h_m(z)]_+^2) = 0.
\end{equation}
Here, the terms are squared to obtain continuously differentiable functions.
In general, the linear independence constraint qualification (LICQ) is not satisfied for such constraints. In our case, for example, the gradient of the rewritten obstacle avoidance constraint is zero at and outside the obstacle boundary. 

Numerous obstacles can be described using the above formulation. For example, any polyhedral set can be cast as a set of affine constraints 
\begin{equation*}
O = \{z \in {\rm I\!R}^{n_z} : b_i - a_i^\intercal z>0, \hspace{0.1cm} i=1,...,m\}.
\end{equation*}
Also non-convex polytopes, such as a cross shaped obstacle, can often be constructed as the union of a set of intersecting convex polygons.
Another type of obstacle shape that can be considered are balls and ellipsoids, given by 
\begin{equation*}
O = \{z \in {\rm I\!R}^{n_z} : 1-(z-c)^\intercal E(z-c)>0, \hspace{0.1cm} i=1,...,m\}.
\end{equation*}
Furthermore, sections of discs can be described as the combination of an outer radius constraint, an inner radius constraint, and a separating hyperplane (affine) constraint. For example, a half-disc obstacle is shown in Figure \ref{fig:IllusHoldAndIntPoint}.

In addition, a set of polynomial functions $h_i(z)$ can be used to define a more general semi-algebraic set. Finally, other functions, such as trigonometric functions, are also possible, as long as they are continuously differentiable.

\subsection{Vehicle models} \label{subsec:Vehicle models}
The problem under consideration is the real-time computation of the optimal trajectory for a motion system. This motion system can be a robot, a satellite, a car, etc. For the remainder of this paper, this system will be called \lq vehicle\rq, as the example models discussed below will be of the vehicle type.

The vehicle is described by a state vector $q$ denoting its position and orientation. In this paper, the considered geometry is two-dimensional. The state therefore has three components: position components $x$ and $y$, and a heading angle $\theta$. The vehicle is steered by control inputs $u$, and its system dynamics are governed by nonlinear ordinary differential equations $\dot{q} = f(q, u)$.

Two vehicle models will be used in the numerical validation of the proposed methodology, cf. Section \ref{sec: Numerical Simulations}. The first vehicle model is the simple bicycle model \citep{rajamani2011vehicle}, where slip of the wheels is neglected. A bicycle is controlled by two control inputs, the velocity $v$ and the steering angle of the front wheel(s) $\delta_f$. The corresponding equations of motion are
\begin{align} \label{eq:Bicycle}
\dot{x} &=  v\cdot \mathrm{cos}(\theta) \nonumber \\ 
\dot{y} &=  v\cdot \mathrm{sin}(\theta) \\ 
\dot{\theta} &=  \frac{v}{L} \mathrm{tan}(\delta_f). \nonumber
\end{align} 
Here, $L$ is the distance between the centers of mass of the wheels of the bicycle.

The second nonlinear vehicle model considered in this paper is a simplified trailer model \citep{embedded}. Again, slip of the wheels is neglected. This model's inputs are the velocity reference $u_x$ and $u_y$ of the towing vehicle. This velocity reference is tracked by a low-level velocity controller. The equations of motion of the trailer model are 
\begin{align}  \label{eq:Trailer}
\dot{x} &=  u_x + L\mathrm{sin}(\theta)\cdot\dot{\theta} \nonumber \\ 
\dot{y} &=  u_y - L\mathrm{cos}(\theta)\cdot\dot{\theta} \\ 
\dot{\theta} &=  \frac{1}{L} (u_y \mathrm{cos}(\theta) - u_x \mathrm{sin}(\theta)). \nonumber
\end{align} 
Here, $L$ is the distance between the center of mass of the trailer vehicle and the fulcrum connecting to the towing vehicle.

\subsection{NMPC formulation} \label{subsec:NMPC formulation}
The continuous-time dynamics describing the motion of the system are discretized using a nonlinear integrator, in this case a fourth order explicit Runge-Kutta method, resulting in the discrete-time representation:
\begin{equation}
q_{k+1} = \varPhi_k(q_k, u_k)
\end{equation}

The NMPC problem is the following:
\begin{align} \label{eq:Original Problem}
{\textrm{\textbf{minimize }}}\hspace{0.2cm}& \ell_N(q_N) + \displaystyle\sum^{N-1}_{k=0}\ell_k(q_k,u_k), \\
{\textrm{\textbf{subject to }}}\hspace{0.2cm}&q_{k+1} = \Phi_k(q_k, u_k), \hspace{0.1cm} k=0,...,N-1, \\
& \psi_{i}(z_k) = 0, \hspace{0.1cm} i=1,...,N_O, \hspace{0.1cm} k=1,...,N, \\
&u_k \in U,  \hspace{0.1cm} k=0,...,N-1,
\end{align}
where $N$ denotes the horizon length and $N_O$ the number of obstacles. 
The obstacle cost functions $\psi_i$ are defined as in (\ref{eq:Obstcost}), and the position $z_k$ is a subvector of the state vector $q_k$. The stage costs are quadratic functions expressing the distance of the state and input variables to the reference state and input: 
\begin{equation*}
\ell_k(q_k,u_k) = (q_k-q_{\mathrm{ref}})^\intercal Q_k(q_k-q_{\mathrm{ref}}) + (u_k-u_{\mathrm{ref}})^\intercal R_k(u_k-u_{\mathrm{ref}}),
\end{equation*}
with the terminal cost 
\begin{equation*}
\ell_N = (q_N-q_{\mathrm{ref}})^\intercal Q_N(q_N-q_{\mathrm{ref}}).
\end{equation*}
The matrices $Q_N$, $Q_k$ and $R_k$ are positive (semi-)definite matrices. A set of input constraints $U$ on which it is easy to project, can straightforwardly be accounted for by the PANOC algorithm. Typical constraints of this type are box constraints of the form  $U = \{u \in {\rm I\!R}^{n_u} : u_\mathrm{min} \leq u \leq u_\mathrm{max}\}$.

\section{METHODOLOGY}

This section discusses the method employed for solving problem (\ref{eq:Problem}), which consists of four parts: (i) a reformulation  of the optimization problem itself; (ii), an optimization algorithm for solving the problem for a fixed value of the penalty parameters; (iii), a penalty method algorithm, which allows for an adequate trade-off between the least-squares objective and the obstacle cost; and (iv), heuristics that facilitate convergence to a trajectory that both reaches the destination and avoids all obstacles.

\subsection{NMPC reformulation} \label{subsec:NMPC reformulation}
Two transformations are applied to the optimization problem before we can introduce a first-order algorithm to solve it. First, the equality constraints representing obstacle avoidance in problem (\ref{eq:Original Problem}) are replaced by appropriate penalty functions in the objective, also known as soft constraints. Given the formulation of the obstacle cost function (\ref{eq:Obstcost}), it is straightforward to construct a quadratic penalty function, $\widetilde{\psi}(z_k) = \frac{1}{2}\mu_k \psi(z_k)^2$, with penalty factors $\mu_k$. This obstacle penalty function has the advantage of being continuously differentiable, in contrast to an exact penalty formulation of these constraints. It is also better conditioned than higher order penalties, with gradient:
\begin{equation}
\nabla \widetilde{\psi}(z_k) = \mu_k\displaystyle\sum^{m}_{i=1} [h_i(z_k)]_+  \displaystyle\prod_{j\neq i} [h_j(z_k)]_+^2 \nabla h_i(z_k).
\end{equation}

Note that \[ \nabla ([w]_+^2) = 2 [w]_+ \nabla ([w]_+) = 2 [w]_+ \nabla w. \] 

Second, the state vectors are eliminated from the optimization problem by integrating the nonlinear kinematic equations of the motion system
\begin{equation}
F_{k+1}(u) = \varPhi_k(F_k(u), u_k), 
\end{equation}
with $F_0(u) = q_0$. This is the so-called single-shooting formulation, where the control inputs are the only remaining decision variables, and the initial state vector is a parameter.

The resulting optimization problem then becomes 
\begin{equation} \label{eq:Problem}
\underset{u \in U^N}{\textrm{\textbf{minimize }}} \ell(u),
\end{equation}
where the objective function is given by
\begin{equation} \label{eq:Objective}
\ell(u) = \ell_N(F_N(u)) + \displaystyle\sum^{N-1}_{k=0}\ell_k(F_k(u),u_k) + \frac{1}{2} \displaystyle\sum^{N}_{k=1} \displaystyle\sum^{N_O}_{i=1} \mu_{ik} {\psi}^2_i(F_k(u)).
\end{equation}
and $U_N=\underbrace{U\times\cdots\times U}_{N\ \textrm{times}}$.

\subsection{PANOC algorithm}

For solving problem (\ref{eq:Problem}) with a fixed value for the penalty parameters $\mu_{ik}$, we employ the recently introduced proximal averaged Newton-type method for optimal control \citep{stella2017simple}. This algorithm, presented in Alg. \ref{alg:PANOC}, achieves a fast convergence by combining proximal gradient and limited memory quasi-Newton (L-BFGS) steps. In this manner, curvature information of the optimization problem is incorporated without calculating second-order derivatives. A set of input constraints $U^N$ can straightforwardly be taken into account via the projection step, step \ref{lst:line:ubar}, in the iterative scheme. 

\begin{algorithm}[H]
\caption{PANOC algorithm for problem (\ref{eq:Problem})}\label{alg:PANOC}
\begin{algorithmic}[1]
\Statex \textbf{Input:} $L_\ell > 0$, $\gamma \in (0,\frac{1}{L_\ell})$, $\sigma \in (0, \frac{\gamma}{2} (1-\gamma\frac{L_\ell}{2}))$, $u^0 \in {\rm I\!R}^{n}$, $\tau > 0$, L-BFGS memory $m$. 
\For{$k = 0,1,2,\ldots$} 
  \State $\bar{u}^k \leftarrow \Pi_{U^N} (u^k - \gamma \nabla \ell(u^k))$  \label{lst:line:ubar}
  \State $r^k \leftarrow \frac{u^k - \bar{u}^k}{\gamma}$  \label{lst:line:r}
  \If {$\|r^k\|_{\infty} < \tau$}  \State {\textbf{stop} with solution $\bar{u}^k$.} 
\EndIf
  \State $d^k = -H_k r^k$ using L-BFGS  
  \State $u^{k+1} \leftarrow u^k - (1 - \alpha_k) \gamma r^k + \alpha_k d^k$, with $\alpha_k$ the largest in $\{ \frac{1}{2^i}: i \in {\rm I\!N} \}$ such that
  \begin{equation}
  \varphi_\gamma(u^{k+1}) \leq \varphi_\gamma(u^k) - \sigma\|r^k\|^2
\end{equation}   
\EndFor 
\end{algorithmic}
\end{algorithm}

In this algorithm, $L_\ell$ denotes the Lipschitz constant of the objective function, $\ell$. If this is not known a priori, as is often the case in practice, the PANOC algorithm can also run with a Lipschitz estimate which is then updated in between iterations, by adding another step after step \ref{lst:line:r}
\begin{algorithmic}[0]
\State {\footnotesize \ref{lst:line:r}bis:} \textbf{if} {$\ell(\bar{u}^k) > \ell(u^k) - \gamma \nabla \ell(u^k)^\intercal r^k + \frac{L_\ell}{2} \|\gamma r^k\|^2$} \textbf{then}
\State \hspace{1.0cm} $\gamma \leftarrow \frac{\gamma}{2}$, $L_\ell \leftarrow 2 L_\ell$, $\sigma \leftarrow \frac{\sigma}{2}$, go to step \ref{lst:line:ubar}.   
\end{algorithmic}
 In addition, two new functions were introduced. The first is the fixed-point residual operator
\begin{equation}
R_\gamma(u) = \frac{1}{\gamma} (u - \Pi_{U^N}(u - \gamma \nabla \ell(u)))
\end{equation}
The second is the forward-backward envelope, first introduced by \citet{patrinos2013proximal}, which can be computed as 
\begin{equation}
\varphi_\gamma(u) = \ell(u) - \frac{\gamma}{2}\|\nabla\ell(u)\|^2 + \mathrm{dist}_U^2(u - \gamma \nabla\ell(u))).
\end{equation} 

For a more in-depth discussion on the properties of PANOC, the reader is referred to \citep{stella2017simple}.
\subsection{Penalty method algorithm}

Given the definition of the obstacle constraint function (\ref{eq:Obstcost}), the obstacles are completely avoided when for every point $z$ of the trajectory and for all obstacles, $\psi_O(z) = 0$ holds. However, solving the optimization problem with the objective as defined in (\ref{eq:Objective}), the solution will likely show a trade-off between low stage costs and low obstacle costs. This trade-off depends on the value of the penalty factors. In order to enforce the obstacle constraints to an acceptable predefined tolerance, we employ a penalty method, as shown in Alg. \ref{alg:Penalty}. Here, it is made explicit that the objective function is parametrized in the penalty factors, hence the notation $\ell(u, \mu)$. Outer iterations are denoted by subscripts, inner iterations by superscripts, and $u$ and $\mu$ represent the vector of control inputs and penalty factors, respectively.

\begin{algorithm}[H]
\caption{Penalty method for problem (\ref{eq:Problem})}\label{alg:Penalty}
\begin{algorithmic}[1]
\Statex \textbf{Input:} $u^0_0 \in {\rm I\!R}^{n}, \mu_0\in{\rm I\!R}^{N_ON}, \eta_* > 0, \tau_*>0, \{\tau_k\} \rightarrow \tau_*, \omega>1 $
\For{$k = 0,1,2,...$} 
  \State Minimize $\ell(u,\mu_k)$ with starting point $u^0_k$, using PANOC with termination criterion $\|R_\gamma (u)\|_{\infty} < \tau_k$ to find $u^*_k$. 
\If {$\|R_\gamma (u^*_k)\|_{\infty} \leq \tau_* \land \|\psi(u^*_k)\| \leq \eta_*$}  \State {\textbf{stop} with solution $u_k^*$.} 
\EndIf
  \State $\mu_{k+1} \leftarrow \omega\mu_k$ 
  \State $u_{k+1}^0 \leftarrow u_k^*$ \label{lst:line:warm starting}
\EndFor 
\end{algorithmic}
\end{algorithm}

In the penalty method, the penalty factors are raised until the optimization solver has converged to a solution for which the norm of the obstacle cost function is lower than a certain tolerance $\eta_*$. The quadratic penalty method is well known to only be exact ($\eta_* = 0$) if the penalty factors are equal to infinity  \citep{nocedal2006numerical}. Therefore, a strictly positive tolerance is chosen. In our algorithm, this is often on the order of $10^{-2}$. %
Virtual enlargements of the obstacle complement this formulation, so that the real obstacles can in fact be completely avoided, even though the constraint tolerance is strictly positive. Such enlargements also allow the formulation to be used for a vehicle with a finite width.

The penalty update factor $\omega$ is used to increase the penalty factors at each outer iteration. In practice, appropriate values of this factor lie between $2-10$ and there is always a trade-off in choosing this value: low values make the different optimization problems more similar and thus easier to warm-start, but more problems will have to be solved in order to converge to a feasible solution. High values in contrast, render the consecutive optimization problems more difficult, but fewer of them are needed. It is observed that for our motion planning problem, the optimization problems do not suffer from a high penalty update factor, thus $\omega$ is here chosen to be 10. In addition, after every update of the penalty factor, the solver is warm-started with the solution from the previous iteration, step \ref{lst:line:warm starting}. After every MPC step, it is common practice to warm-start the next optimal control problem, shifting the vector of control inputs over one time instant and adding an initial guess, often the zero vector, for the last time instant. Similarly, the penalty factors are shifted, and a vector of ones is added for the last time instant.

\begin{figure}
  \centering
  \begin{subfigure}[b]{0.49\columnwidth}
    \includegraphics[width=\textwidth]{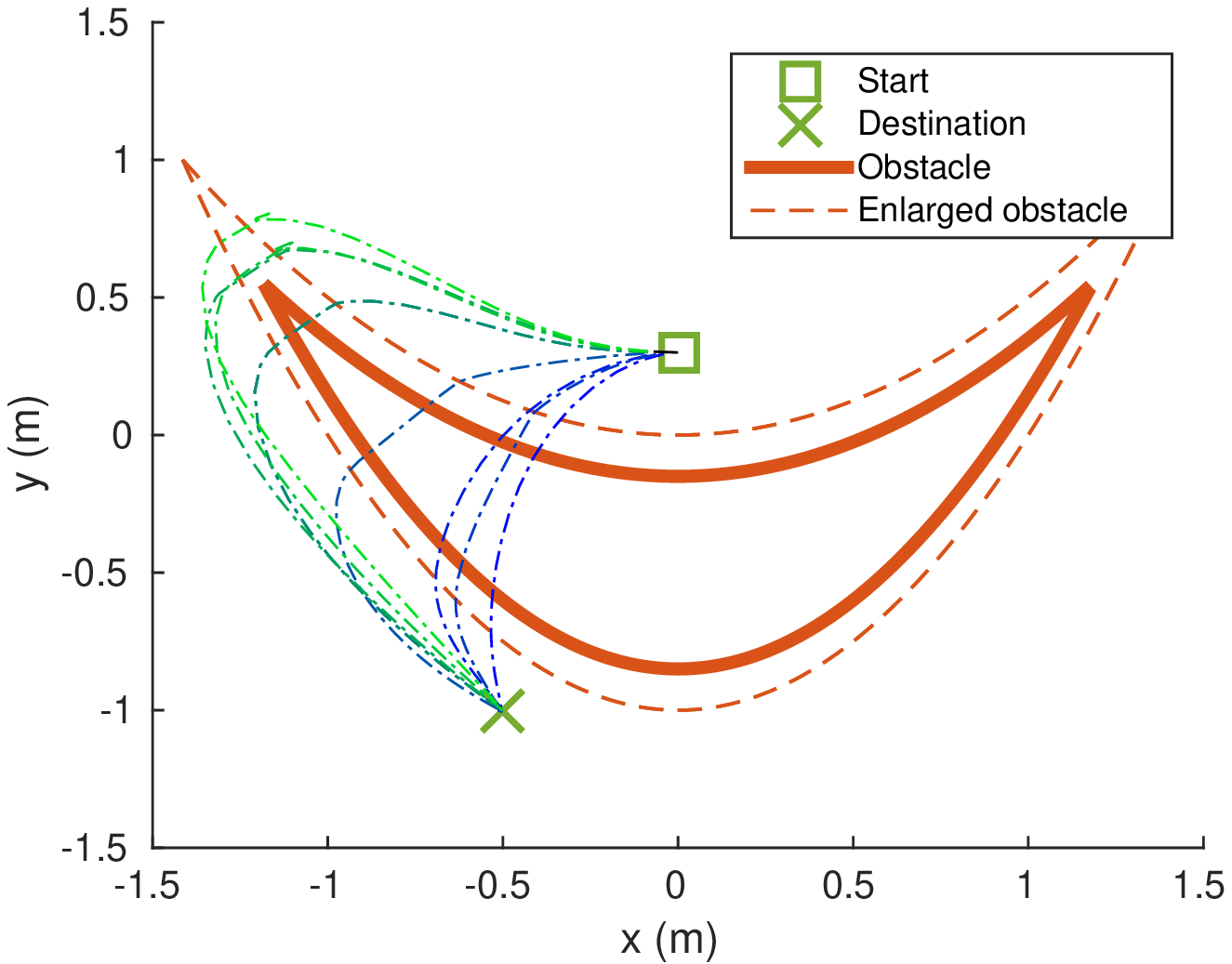}
	\caption{Penalty method.}
    \label{fig:IllusPenaltyMethod}
  \end{subfigure}
  \begin{subfigure}[b]{0.49\columnwidth}
    \includegraphics[width=\textwidth]{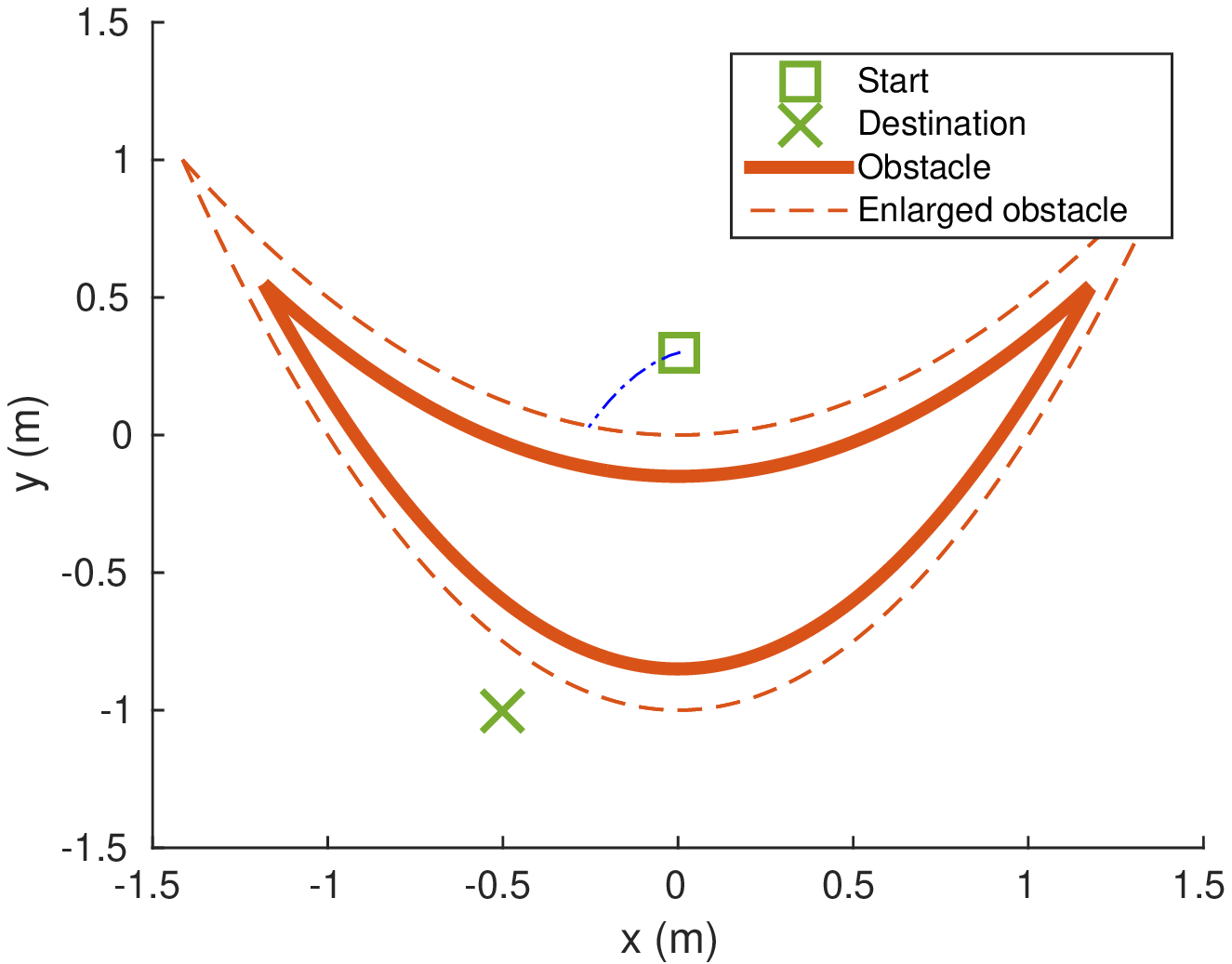}
    \caption{Fixed (high) penalty factor.}
    \label{fig:IllusPenaltyMethod2}
  \end{subfigure}
\caption{\small Illustration of the penalty method. The enlarged obstacle is defined by $O = \{(x,y):y>x^2, y< 1+x^2/2\}$.} 
\label{fig:IllusPenaltyMethod1and2}
\end{figure}

An illustration of the penalty method applied to a problem with a crescent obstacle is shown in Figure \ref{fig:IllusPenaltyMethod}. The trajectories ranging from blue to green correspond to increasing penalty factors, which determine the balance between a feasible trajectory and one that is optimal for the least squares objective. The enlarged obstacle can never be completely avoided, but for high enough penalty factors, the original obstacle is avoided, as illustrated by the final two green trajectories. The combination  of a virtual enlargement of the obstacle and finite values for the penalty factors is therefore indeed successful. Figure \ref{fig:IllusPenaltyMethod} also demonstrates that successive trajectories are usually similar in shape even though the penalty factors are updated somewhat aggressively in this paper. Therefore, warm-starting aids tremendously in the convergence of the penalty method. 

The application of the penalty method may have an additional benefit for obstacle avoidance problems, also illustrated by Figure \ref{fig:IllusPenaltyMethod}. Assuming some obstacle is blocking the shortest path from start to destination, the initial trajectory calculated with low penalty factors is very likely to arrive at the destination while violating obstacle avoidance constraints. Subsequent iterations with higher and higher penalty factors tend to push the trajectory to the edge of the obstacle, while remaining connected to the destination. In contrast, solving the problem only once with a high value for the penalty factors can impede convergence to a trajectory that reaches the destination, as the vehicle is more likely to get stuck behind an obstacle, as shown in Figure \ref{fig:IllusPenaltyMethod2}. Using the penalty method for the problems considered in this paper therefore aids in avoiding the local minimum.

\subsection{Heuristics}
\begin{figure}
  \centering
  \begin{subfigure}[b]{0.49\columnwidth}
    \includegraphics[width=\textwidth]{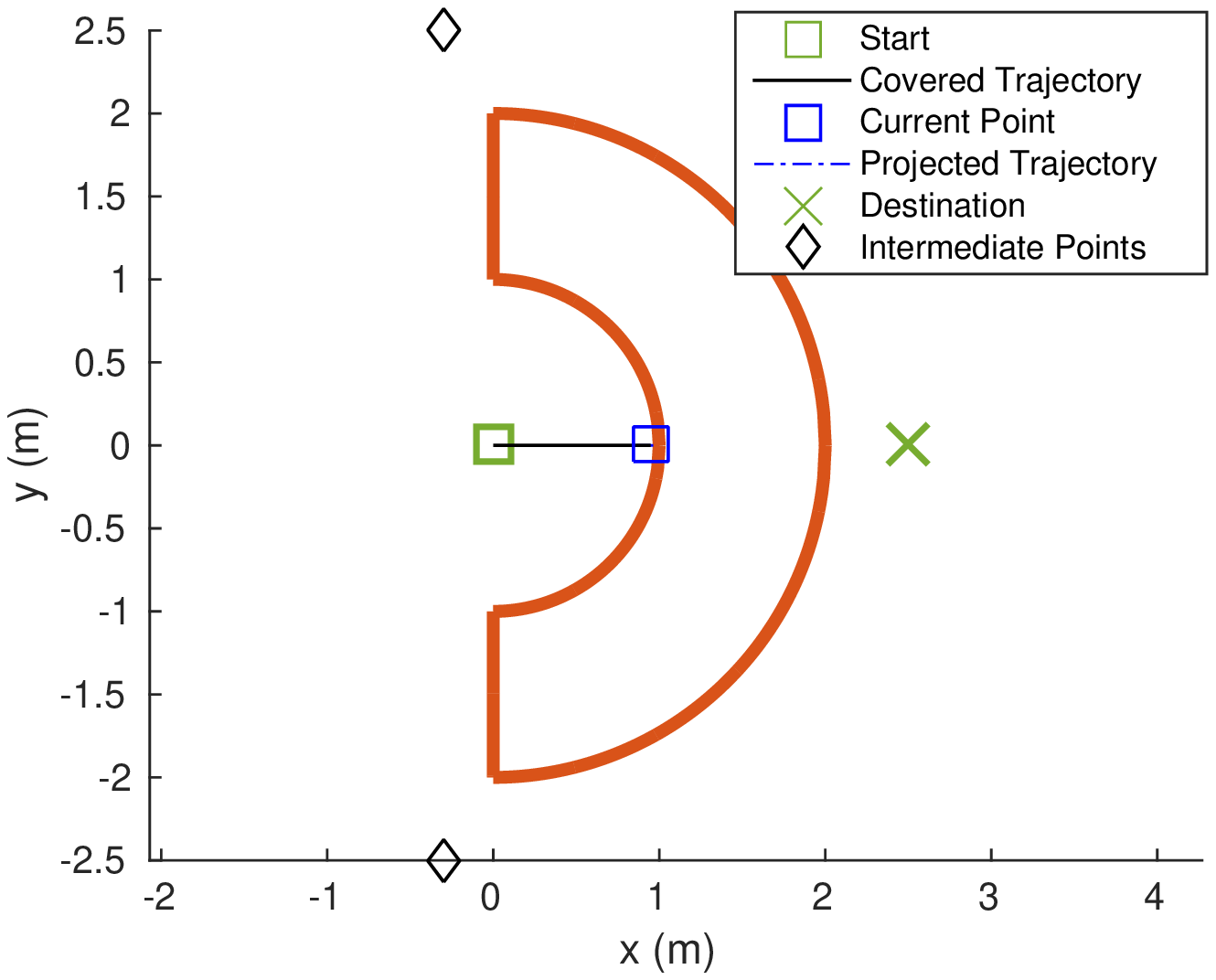}
	\caption{High penalty factors cap.}
    \label{fig:local minimum}
  \end{subfigure}
  \begin{subfigure}[b]{0.49\columnwidth}
    \includegraphics[width=\textwidth]{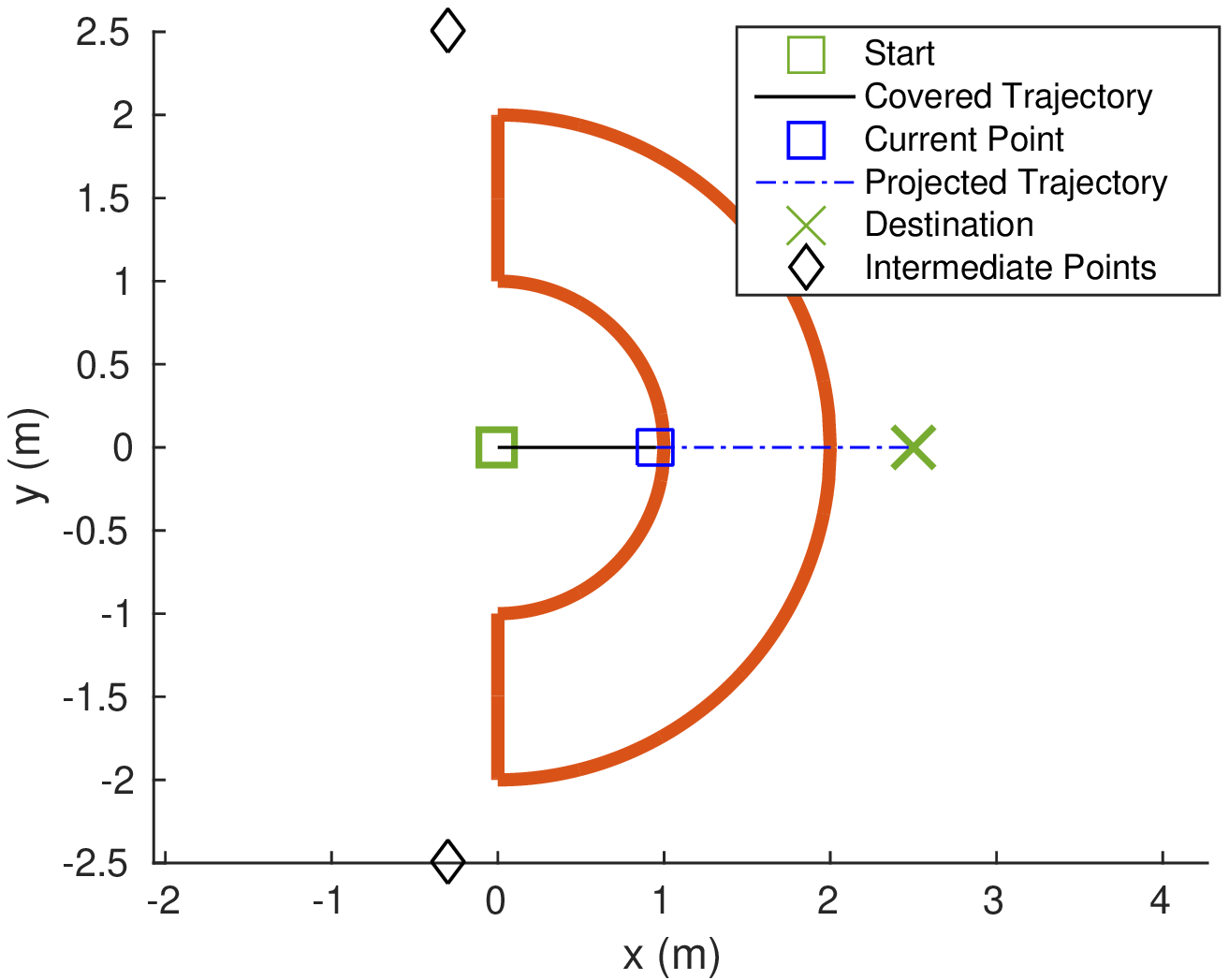}
    \caption{Low penalty factors cap.}
    \label{fig:hold}
  \end{subfigure}
\caption{\small Illustration of the local minimum behind the obstacle, and the hold-in-place heuristic and choice of intermediate points. The half-disc shaped obstacle is defined by $O = \{(x,y): x^2+y^2 > 1, x^2+y^2 < 4, x>0\}$. } 
\label{fig:IllusHoldAndIntPoint}
\end{figure}

Optimization problem (\ref{eq:Problem}) is a general nonlinear, non-convex problem. The obstacles render the solution space non-convex, and thus often local minima exist near obstacles. An example hereof is illustrated in Figure \ref{fig:local minimum}. To aid in the convergence to a feasible trajectory that reaches the destination, three additional heuristics are applied to the algorithm: (i) the penalty factors are capped at an appropriate value; (ii) the vehicle is stopped if the obstacle costs do not satisfy the specified tolerance  $\eta_*$; and (iii) whenever the vehicle remains in place for more than one time instant, it is guided to an intermediate destination before continuing on towards the final destination. Below, the rationale for each of these heuristics is explained.

Large penalty factors render the problem ill-conditioned, which impedes fast convergence of the method. In particular, PANOC uses an estimate of the Lipschitz constant of the objective function to determine the stepsize, and the Lipschitz constants of the obstacle cost penalties scale linearly with the penalty factors. Therefore, these factors are capped at a reasonably low value, for example $10^4$. With this cap, however, the algorithm is not guaranteed to find a solution for which the obstacle costs are sufficiently low. 

In order to solve this problem, the following heuristic is applied: If the obstacle cost is not smaller than the prescribed tolerance within the next three time steps, the vehicle performs an emergency stop. This strategy prevents collisions with obstacles, but causes the vehicle to be more likely to get stuck. Figure \ref{fig:hold} illustrates this strategy in case of a half-disc shaped obstacle. Starting from the green square, the vehicle moves to the blue square in seven MPC steps. There, the calculated trajectory threatens to violate the obstacle constraint within the next three time steps, because the corresponding penalty factors are too low. Instead of following this trajectory, the vehicle is stopped.

Finally, to assist the vehicle in circumventing all obstacles, another heuristic is implemented. If the vehicle is stuck behind an obstacle, then the reference state is temporarily replaced by an intermediate destination. The intention behind this is to guide the vehicle around the obstacle. A good intermediate destination is easy to reach from both the point where the vehicle was previously stuck and the final destination. It will usually be close to a corner or edge point of the obstacle. This principle is also illustrated in Figure \ref{fig:IllusHoldAndIntPoint}.  Appropriate intermediate destinations lie near the black diamonds.

To avoid getting stuck in a local optimum near an obstacle, a suitable set of intermediate destinations must be available and an appropriate choice from this set of points is necessary. The user may provide such a set, based on knowledge of the obstacle definitions (and therefore the locations of their corners). This approach, however, is hard to justify in an automated setup. 
In order to generate suitable intermediate points automatically, variants of Dijkstra's algorithm can be used to perform a simple graph search. This paper utilizes the A* search algorithm \citep{hart1968formal}, because of its simplicity and efficiency. The worst-case complexity of this algorithm in case a consistent heuristic cost is used, is $O(N)$ \citep{martelli1977complexity}, with $N$ the number of nodes in the graph. The heurstic cost used here is the Euclidean distance between a node and the goal node, which is indeed consistent. 

\begin{figure}
\begin{center}
\includegraphics[width=5.4cm]{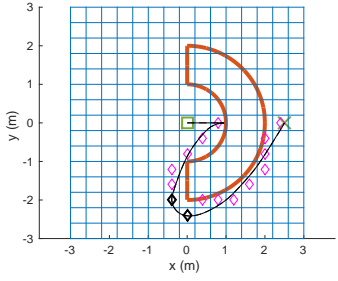}    
\caption{\small Illustration of the graph search, represented by the magenta triangles. At both of the black triangles the direction changes (first left-right, then up-down), so these comprise the set of intermediate destinations. The black line denotes the final trajectory.} 
\label{fig:GraphSearchIllustration}
\end{center}
\end{figure}

Figure \ref{fig:GraphSearchIllustration} shows that the graph search returns a feasible trajectory from the current point to the destination. Intermediate points can be extracted from this trajectory by moving through it and recording points at which the left-right or up-down direction switches. In an unobstructed space, a graph search would find a straight path. Hence, direction changes stem from the presence of an obstacle, and the points at which they occur are close to the corners of the obstacle. These points are therefore suitable candidates for intermediate destinations. In the simulations below, this approach is used and each point of the set of intermediate destinations is successively visited. When this set is exhausted, the reference state is reset to original destination.

\section{NUMERICAL SIMULATIONS} \label{sec: Numerical Simulations}

The proposed methodology is illustrated and analyzed by means of numerical simulations for a wide variety of obstacle configurations, and two different vehicle models. All simulations were performed on a notebook with Intel(R) Core(TM) i7-7600U CPU @ 2.80GHz x 2 processor and 16 GB of memory.

\begin{figure}
  \centering
  \begin{subfigure}[t]{0.47\columnwidth}
    \includegraphics[width=\textwidth]{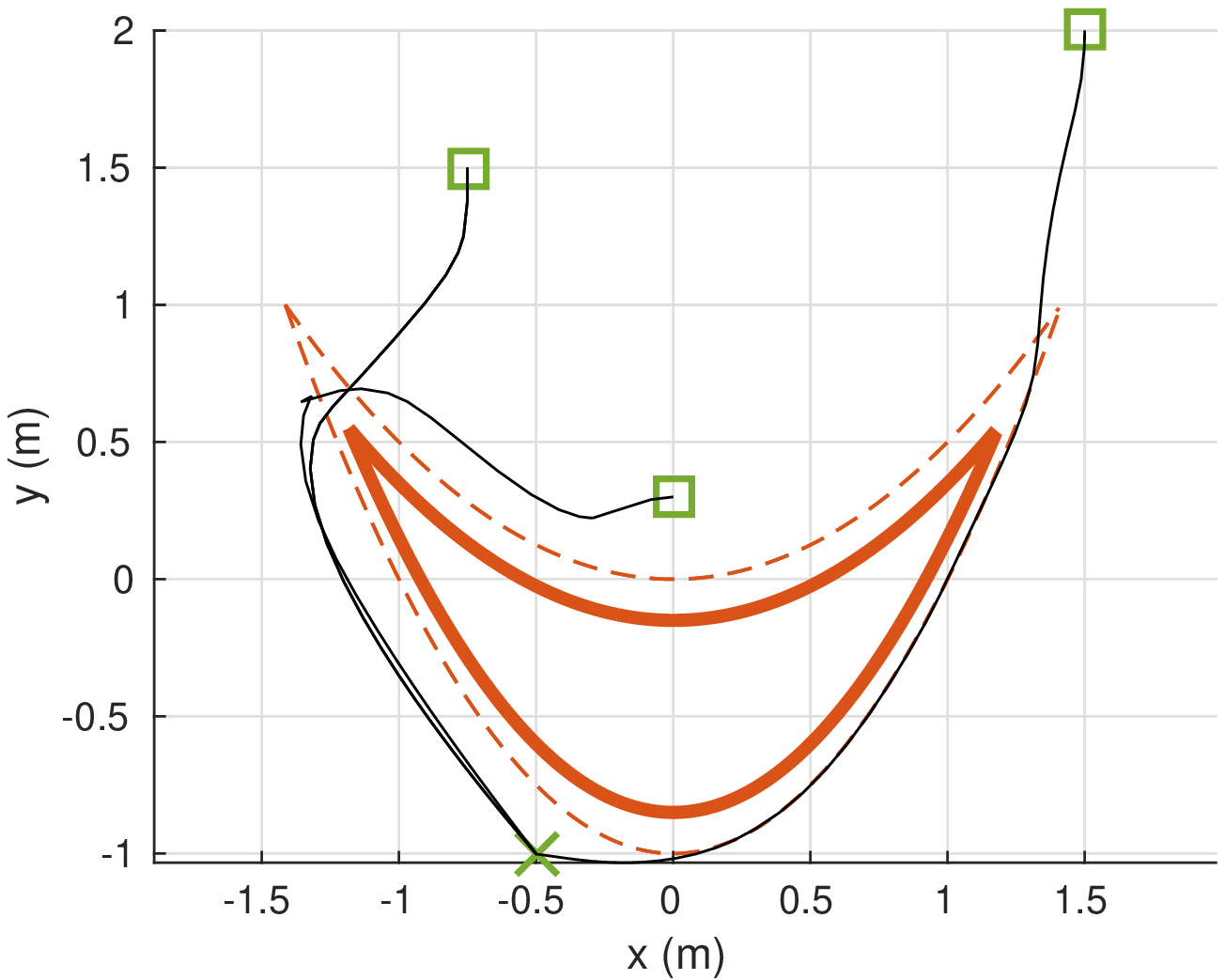}
    \caption{Enlarged obstacle defined as $\{(x,y):y>x^2, y< 1+x^2/2\}$.}
    \label{fig:Crescent}
  \end{subfigure} \quad
    \begin{subfigure}[t]{0.47\columnwidth}
    \includegraphics[width=\textwidth]{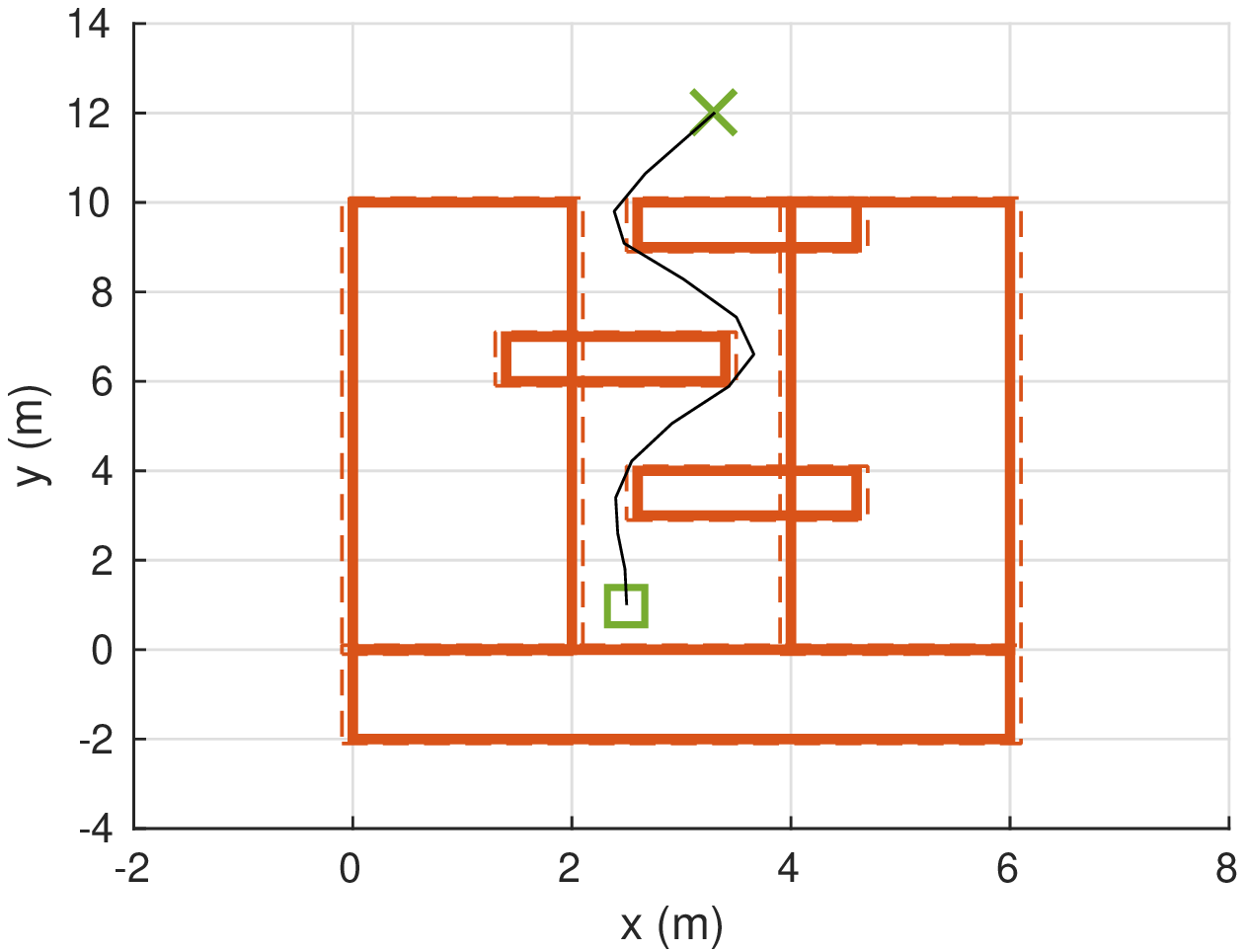}
    \caption{Slightly complex obstacle configuration using rectangular obstacles.}
    \label{fig:Labyrinth}
  \end{subfigure}
  \caption{\small Two obstacle configurations using the trailer model. The conventions for start, destination and obstacles are the same as those in Figure \ref{fig:IllusPenaltyMethod}, and the black lines denote the trajectories that were found in each case.}
  \label{Fig:simulationsTrailer}
\end{figure}

\begin{figure}
\centering
  \begin{subfigure}[t]{.47\columnwidth}
  \centering
    \includegraphics[width=\linewidth]{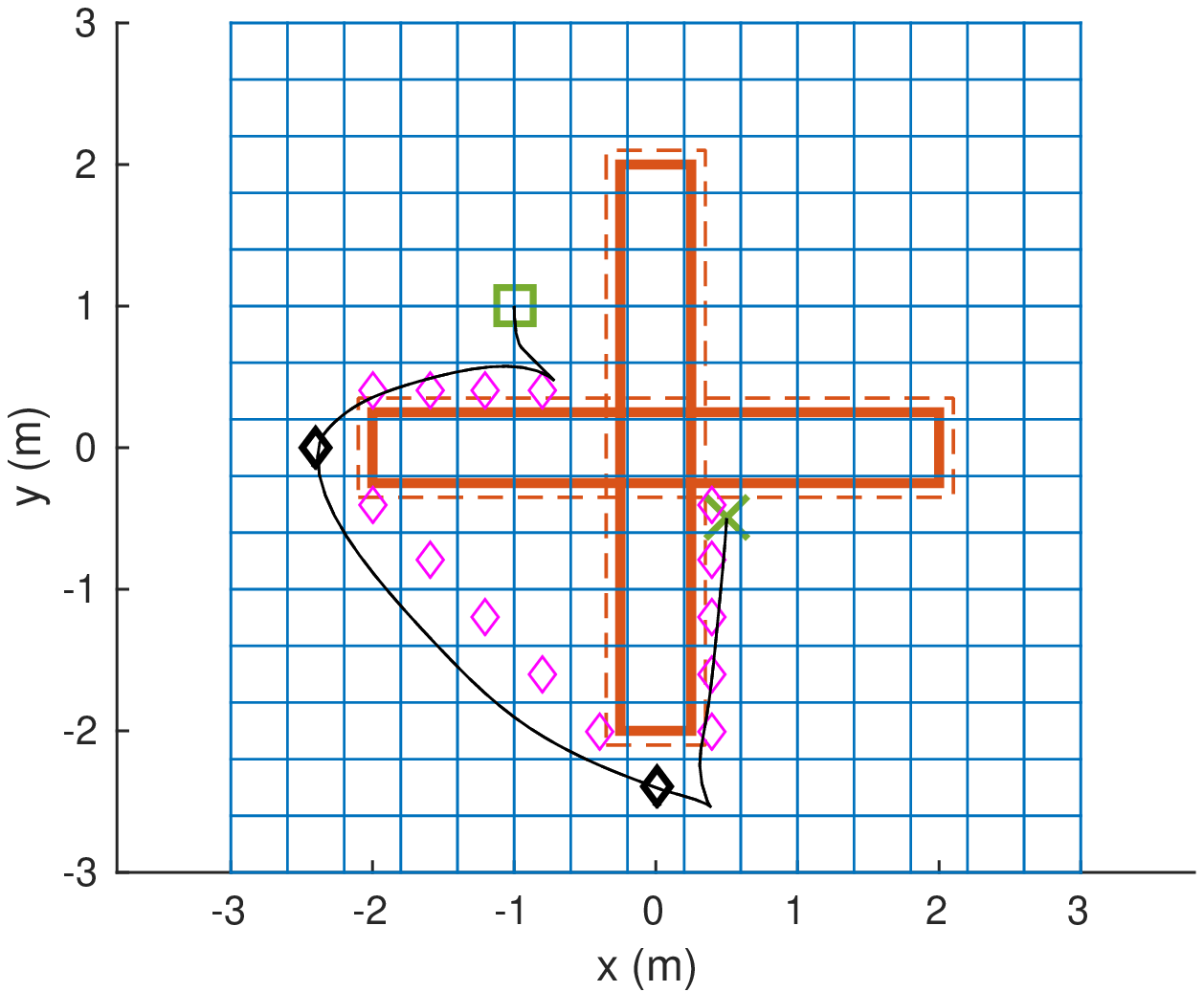}
    \caption{Cross-shaped obstacle as the combination of two rectangles.}
  \end{subfigure} \quad
  \begin{subfigure}[t]{.47\columnwidth}
  \centering
    \includegraphics[width=\linewidth]{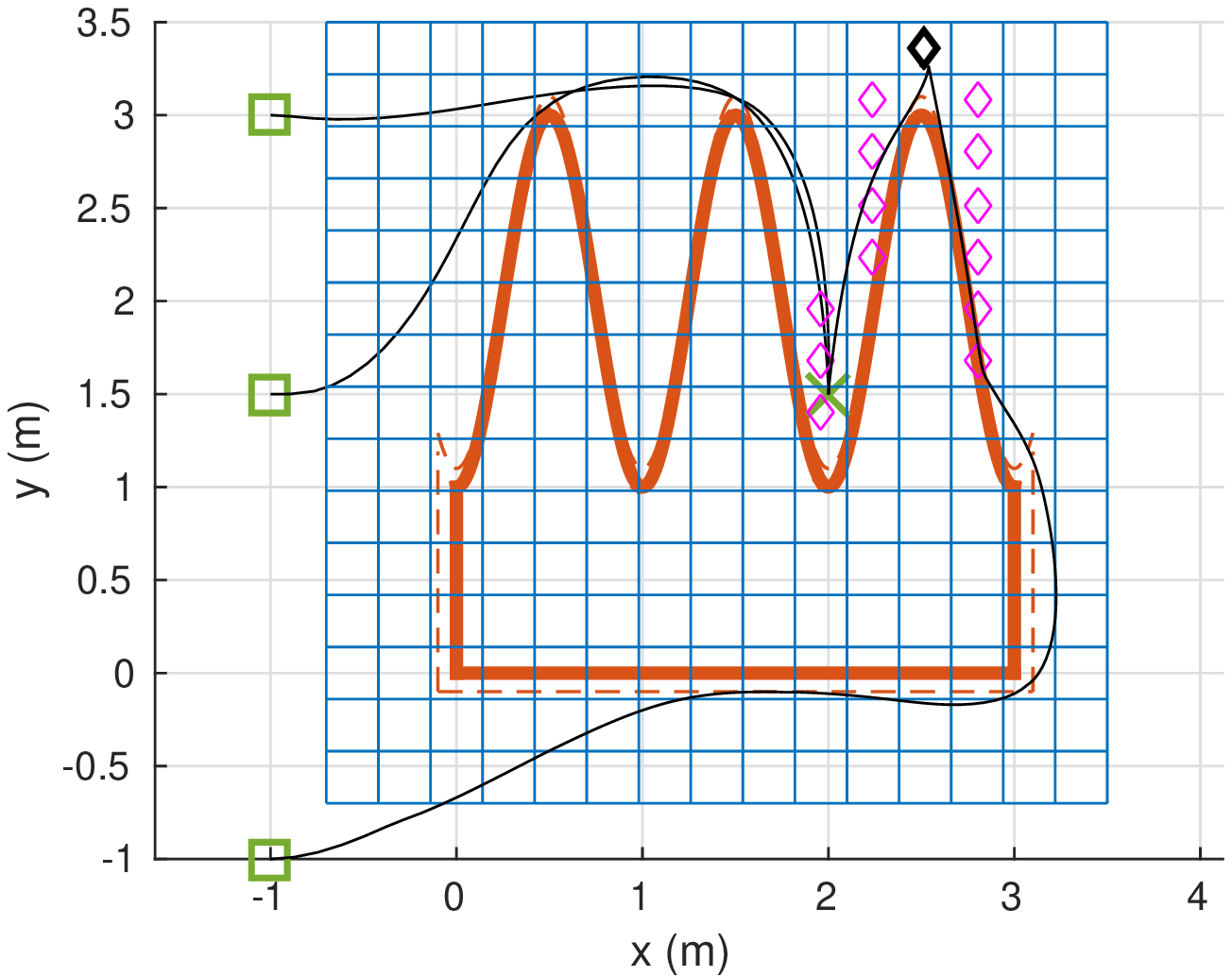}
    \caption{Rack-shaped obstacle, defined by $O = \{(x,y): y < \mathrm{sin}(2\pi x - \pi/2) + 2, y > 0, 0 < x < 3 \}$.}
  \end{subfigure}
  \caption{\small Illustration of the graph search heuristic for two obstacle configurations using the trailer model.}
  \label{Fig:simulationsTrailerGraphSearch}
\end{figure}

Figure \ref{Fig:simulationsTrailer} displays the first set of obstacle configurations, overcome by a vehicle with a trailer. The kinematics of the trailer model, given by (\ref{eq:Trailer}), are discretized using an explicit fourth order Runge Kutta method. The distance between the center of mass of the trailer and the fulcrum of the towing vehicle is $L = 0.5$m. The optimal control problems are solved with sampling time $t_s = 30$ms for Fig. \ref{fig:Crescent} and $t_s = 200$ms for Fig. \ref{fig:Labyrinth}, and horizon length $N = 50$. The inputs are constrained by box constraints 
$-4\mathrm{m/s} \leq u_x,u_y \leq 4\mathrm{m/s}$ at every time instant.
The penalty factors are capped off at $10^4$. It is usually difficult for a first order method, such as PANOC, to find a solution for a strict tolerance, say $10^{-6}$. Therefore, the tolerance $\tau_*$ is set to $10^{-3}$. We observed in simulations that the closed-loop performance is not impacted by this choice, and that a stricter tolerance would be unnecessary. Figure \ref{Fig:simulationsTrailerGraphSearch} displays two more obstacle configurations for which the graph search heuristic proved necessary to find a trajectory that reached the destination.

\begin{figure}
\centering
\begin{subfigure}{.47\columnwidth}
  \centering
  \includegraphics[width=\linewidth]{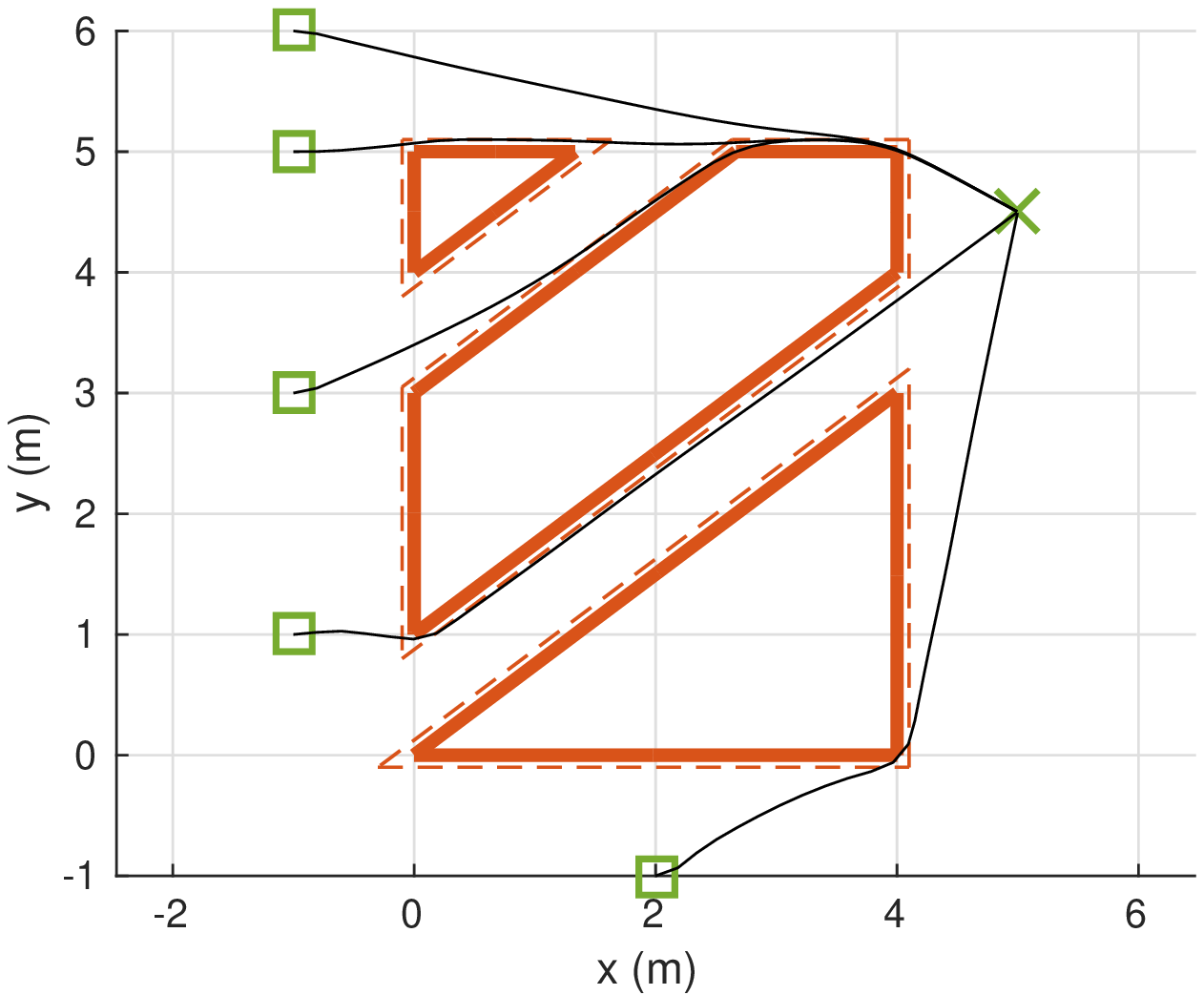}
  \caption{Three polyhedral obstacles that form corridors.}
  \label{fig:Corridors}
\end{subfigure} \quad
\begin{subfigure}{.47\columnwidth}
  \centering
  \includegraphics[width=\linewidth]{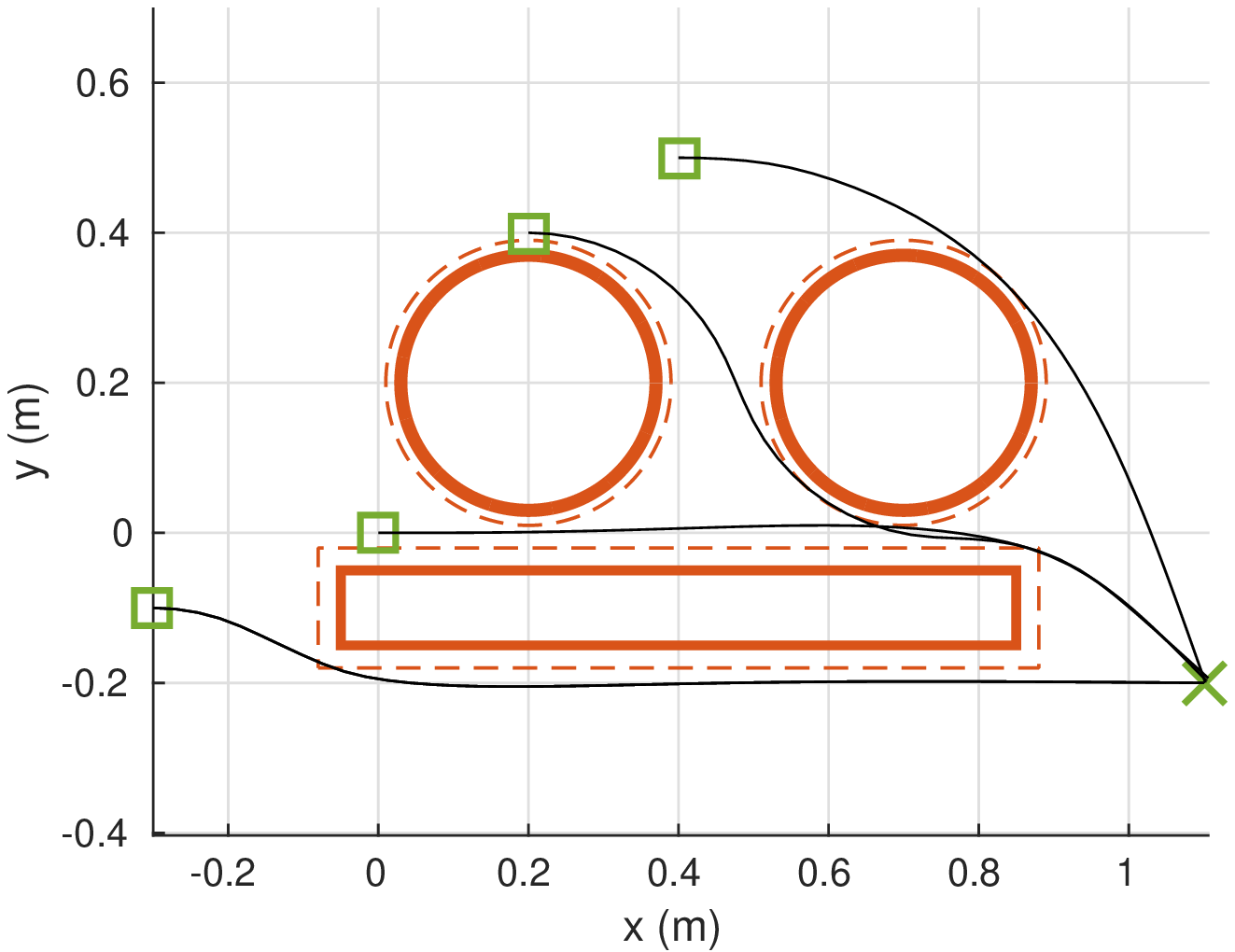}
  \caption{Configuration of one rectangular and two circular obstacles.}
  \label{fig:TwoCircOneRect}
\end{subfigure}
\caption{\small Two obstacle configurations using the bicycle model.}
\label{Fig:simulationsBicycle}
\end{figure}

Similarly, Figure \ref{Fig:simulationsBicycle} displays two obstacle configurations overcome by a vehicle modeled as a bicycle. The relevant kinematics are given by (\ref{eq:Bicycle}). The optimal control problem for these simulations is constructed with the same parameters as before, with the following exceptions: the sampling time $t_s = 50$ms, and the input constraints are $-0.1 \mathrm{m/s} \leq v \leq 4 \mathrm{m/s}$ and $-\frac{\pi}{3} \leq \delta_f \leq \frac{\pi}{3}$ at every time instant. The figures in this section show the versatility of the proposed approach for constructing and solving collision-avoidance problems.

\begin{figure}
\begin{center}
\includegraphics[width=6.4cm]{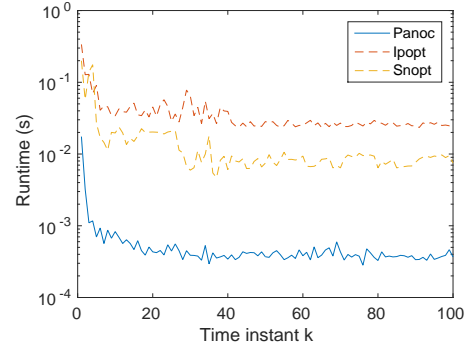}    
\caption{\small Runtime comparison of three solvers: PANOC, SNOPT and IPOPT. This comparison is for the problem of Fig. \ref{fig:TwoCircOneRect}, with initial point $x_0 = (0, 0, 0)$. All solvers were warm-started with their previous solution after each MPC step. The tolerance for all solvers was set to $10^{-3}$.} 
\label{fig:RuntimeComparison}
\end{center}
\end{figure}

Figure \ref{fig:RuntimeComparison} compares the runtime of the proposed methodology with that of state-of-the-art SQP \citep[SNOPT]{gill2005snopt} and IP \citep[IPOPT]{wachter2006implementation} solvers. In these solvers, the obstacle avoidance constraint is incorporated as an inequality constraint, $\psi^2(z) \leq \eta_*^2$. The penalty method algorithm using PANOC clearly outperforms the other solvers, by approximately two orders of magnitude.

\begin{table}[!ht]
\begin{center}
\caption{\small Comparison of closed-loop costs.}\label{tb:closed-loop cost}
\begin{tabular}{|l|c|c|c|}
  \hline
\backslashbox{Example}{Solver}   & PANOC & IPOPT & SNOPT \\ \hline \hline
  Figure \ref{fig:TwoCircOneRect}, $x_0 = (0,0,0)$ & 3.56 & 4.02 & 3.54 \\ \hline
  Figure \ref{fig:Crescent}, $x_0 = (0, 0.3, \pi)$ & 20.78 & 49.43 & 47.50\\ \hline
  Figure \ref{fig:Corridors}, $x_0 = (-1, 1, 0)$ & 17.65 & 19.62 & 26.17\\
  \hline 
    Figure \ref{fig:Corridors}, $x_0 = (-1, 3, 0)$ & 13.21 & 12.76 & 33.48\\
  \hline 
\end{tabular}
\end{center}
\end{table}

The proposed approach is not only faster, but also more adept at finding high quality solutions than the other state-of-the-art solvers. This is illustrated in Table \ref{tb:closed-loop cost}, which lists the closed-loop costs for different problem scenarios. Sometimes, the solvers all converged to the same trajectory, such as for Fig. \ref{fig:TwoCircOneRect}. This was, however, not always the case. For example, IPOPT and SNOPT were both temporarily stuck in the local minimum behind the crescent-shaped obstacle, Fig. \ref{fig:Crescent}, whereas our approach found the optimal trajectory in the first optimal control problem. In other cases, such as the corridors example, Fig. \ref{fig:Corridors}, different paths were taken. Table \ref{tb:closed-loop cost} shows that for these cases, the proposed methodology found trajectories with closed-loop costs as good as or better than the other solvers.

\section{Conclusion}

This paper presents a penalty method framework for solving optimal control problems with collision-avoidance constraints that typically arise in motion planning problems. The application of the penalty method, coupled with virtual enlargements, allows for the avoidance of obstacles of complex geometry. It also benefits the convergence to a trajectory that reaches the destination, by gradually finding a path around the obstacles as the penalty factors are successively increased. In addition, several heuristics are employed, which have been observed to improve convergence to a feasible trajectory that reaches the destination. The resulting optimization problems are solved with PANOC, a first-order method which exhibits low runtime. Numerical simulations with nonlinear vehicle dynamics show the versatility of the proposed approach in solving motion planning problems with general obstacle avoidance constraints. In the limited number of cases considered, the proposed algorithm outperforms state-of-the-art SQP and IP solvers in both runtime and robustness.


\small
\bibliography{ifacconf}            

\begin{thebibliography}{24}
\providecommand{\natexlab}[1]{#1}
\providecommand{\url}[1]{\texttt{#1}}
\providecommand{\urlprefix}{URL }
\expandafter\ifx\csname urlstyle\endcsname\relax
  \providecommand{\doi}[1]{doi:\discretionary{}{}{}#1}\else
  \providecommand{\doi}{doi:\discretionary{}{}{}\begingroup
  \urlstyle{rm}\Url}\fi

\bibitem[{Andersson et~al.(In Press, 2018)Andersson, Gillis, Horn, Rawlings,
  and Diehl}]{Andersson2018}
Andersson, J.A.E., Gillis, J., Horn, G., Rawlings, J.B., and Diehl, M. (In
  Press, 2018).
\newblock {CasADi} -- {A} software framework for nonlinear optimization and
  optimal control.
\newblock \emph{Mathematical Programming Computation}.

\bibitem[{Boyd and Vandenberghe(2004)}]{boyd2004convex}
Boyd, S. and Vandenberghe, L. (2004).
\newblock \emph{Convex optimization}.
\newblock Cambridge university press.

\bibitem[{Debrouwere et~al.(2013)Debrouwere, Van~Loock, Pipeleers, Diehl,
  De~Schutter, and Swevers}]{debrouwere2013time}
Debrouwere, F., Van~Loock, W., Pipeleers, G., Diehl, M., De~Schutter, J., and
  Swevers, J. (2013).
\newblock Time-optimal path following for robots with object collision
  avoidance using lagrangian duality.
\newblock In \emph{Robot Motion and Control (RoMoCo), 2013 9th Workshop on},
  186--191. IEEE.

\bibitem[{Ge and Cui(2002)}]{ge2002dynamic}
Ge, S.S. and Cui, Y.J. (2002).
\newblock Dynamic motion planning for mobile robots using potential field
  method.
\newblock \emph{Autonomous robots}, 13(3), 207--222.

\bibitem[{Gill et~al.(2005)Gill, Murray, and Saunders}]{gill2005snopt}
Gill, P.E., Murray, W., and Saunders, M.A. (2005).
\newblock {SNOPT}: An {SQP} algorithm for large-scale constrained optimization.
\newblock \emph{SIAM review}, 47(1), 99--131.

\bibitem[{Guy et~al.(2009)Guy, Chhugani, Kim, Satish, Lin, Manocha, and
  Dubey}]{guy2009clearpath}
Guy, S.J., Chhugani, J., Kim, C., Satish, N., Lin, M., Manocha, D., and Dubey,
  P. (2009).
\newblock Clearpath: highly parallel collision avoidance for multi-agent
  simulation.
\newblock In \emph{Proceedings of the 2009 ACM SIGGRAPH/Eurographics Symposium
  on Computer Animation}, 177--187. ACM.

\bibitem[{Hart et~al.(1968)Hart, Nilsson, and Raphael}]{hart1968formal}
Hart, P.E., Nilsson, N.J., and Raphael, B. (1968).
\newblock A formal basis for the heuristic determination of minimum cost paths.
\newblock \emph{IEEE transactions on Systems Science and Cybernetics}, 4(2),
  100--107.

\bibitem[{Jerez et~al.(2014)Jerez, Goulart, Richter, Constantinides, Kerrigan,
  and Morari}]{jerez2014embedded}
Jerez, J.L., Goulart, P.J., Richter, S., Constantinides, G.A., Kerrigan, E.C.,
  and Morari, M. (2014).
\newblock Embedded online optimization for model predictive control at
  megahertz rates.
\newblock \emph{IEEE Transactions on Automatic Control}, 59(12), 3238--3251.

\bibitem[{Martelli(1977)}]{martelli1977complexity}
Martelli, A. (1977).
\newblock On the complexity of admissible search algorithms.
\newblock \emph{Artificial Intelligence}, 8(1), 1--13.

\bibitem[{Mercy et~al.(2017)Mercy, Van~Parys, and Pipeleers}]{mercy2017spline}
Mercy, T., Van~Parys, R., and Pipeleers, G. (2017).
\newblock Spline-based motion planning for autonomous guided vehicles in a
  dynamic environment.
\newblock \emph{IEEE Transactions on Control Systems Technology}.

\bibitem[{Montiel et~al.(2015)Montiel, Orozco-Rosas, and
  Sep{\'u}lveda}]{montiel2015path}
Montiel, O., Orozco-Rosas, U., and Sep{\'u}lveda, R. (2015).
\newblock Path planning for mobile robots using bacterial potential field for
  avoiding static and dynamic obstacles.
\newblock \emph{Expert Systems with Applications}, 42(12), 5177--5191.

\bibitem[{Nocedal and Wright(2006)}]{nocedal2006numerical}
Nocedal, J. and Wright, S.J. (2006).
\newblock Numerical optimization 2nd.

\bibitem[{Patrinos and Bemporad(2013)}]{patrinos2013proximal}
Patrinos, P. and Bemporad, A. (2013).
\newblock Proximal newton methods for convex composite optimization.
\newblock In \emph{Decision and Control (CDC), 2013 IEEE 52nd Annual Conference
  on}, 2358--2363. IEEE.

\bibitem[{Patrinos and Bemporad(2014)}]{patrinos2014accelerated}
Patrinos, P. and Bemporad, A. (2014).
\newblock An accelerated dual gradient-projection algorithm for embedded linear
  model predictive control.
\newblock \emph{IEEE Transactions on Automatic Control}, 59(1), 18--33.

\bibitem[{Rajamani(2011)}]{rajamani2011vehicle}
Rajamani, R. (2011).
\newblock \emph{Vehicle dynamics and control}.
\newblock Springer Science \& Business Media.

\bibitem[{Rawlings and Mayne(2009)}]{rawlings2009model}
Rawlings, J. and Mayne, D. (2009).
\newblock \emph{Model Predictive Control: Theory and Design}.
\newblock Nob Hill Pub.

\bibitem[{Richter et~al.(2012)Richter, Jones, and
  Morari}]{richter2012computational}
Richter, S., Jones, C.N., and Morari, M. (2012).
\newblock Computational complexity certification for real-time mpc with input
  constraints based on the fast gradient method.
\newblock \emph{IEEE Transactions on Automatic Control}, 57(6), 1391--1403.

\bibitem[{Sathya et~al.(2018)Sathya, Sopasakis, Van~Parys, Themelis, Pipeleers,
  and Patrinos}]{embedded}
Sathya, A.S., Sopasakis, P., Van~Parys, R., Themelis, A., Pipeleers, G., and
  Patrinos, P. (2018).
\newblock Embedded nonlinear model predictive control for obstacle avoidance
  using panoc.
\newblock In \emph{Proceedings of the 2018 European Control Conference}.

\bibitem[{Scheel and Scholtes(2000)}]{scheel2000mathematical}
Scheel, H. and Scholtes, S. (2000).
\newblock Mathematical programs with complementarity constraints: Stationarity,
  optimality, and sensitivity.
\newblock \emph{Mathematics of Operations Research}, 25(1), 1--22.

\bibitem[{Stella et~al.(2017)Stella, Themelis, Sopasakis, and
  Patrinos}]{stella2017simple}
Stella, L., Themelis, A., Sopasakis, P., and Patrinos, P. (2017).
\newblock A simple and efficient algorithm for nonlinear model predictive
  control.
\newblock In \emph{56th IEEE Conference on Decision and Control}, 1939--1944.

\bibitem[{Stentz(1994)}]{stentz1994optimal}
Stentz, A. (1994).
\newblock Optimal and efficient path planning for partially-known environments.
\newblock In \emph{ICRA}, volume~94, 3310--3317.

\bibitem[{Takahashi and Schilling(1989)}]{takahashi1989motion}
Takahashi, O. and Schilling, R.J. (1989).
\newblock Motion planning in a plane using generalized voronoi diagrams.
\newblock \emph{IEEE Transactions on robotics and automation}, 5(2), 143--150.

\bibitem[{W{\"a}chter and Biegler(2006)}]{wachter2006implementation}
W{\"a}chter, A. and Biegler, L.T. (2006).
\newblock On the implementation of an interior-point filter line-search
  algorithm for large-scale nonlinear programming.
\newblock \emph{Mathematical programming}, 106(1), 25--57.

\bibitem[{Wang and Ding(2014)}]{wang2014synthesis}
Wang, P. and Ding, B. (2014).
\newblock A synthesis approach of distributed model predictive control for
  homogeneous multi-agent system with collision avoidance.
\newblock \emph{International Journal of Control}, 87(1), 52--63.

\end{thebibliography}

\end{document}